# STEIN ESTIMATION FOR THE DRIFT OF GAUSSIAN PROCESSES USING THE MALLIAVIN CALCULUS

BY NICOLAS PRIVAULT AND ANTHONY RÉVEILLAC

*City University of Hong Kong and Université de la Rochelle*

We consider the nonparametric functional estimation of the drift of a Gaussian process via minimax and Bayes estimators. In this context, we construct superefficient estimators of Stein type for such drifts using the Malliavin integration by parts formula and superharmonic functionals on Gaussian space. Our results are illustrated by numerical simulations and extend the construction of James–Stein type estimators for Gaussian processes by Berger and Wolpert [*J. Multivariate Anal.* **13** (1983) 401–424].

**1. Introduction.** The maximum likelihood estimator $\hat{\mu}$ of the mean $\mu \in \mathbb{R}^d$ of a Gaussian random vector $X$ in $\mathbb{R}^d$ with covariance $\sigma^2 I_{\mathbb{R}^d}$ under a probability $\mathbb{P}_\mu$ is well known to be equal to $X$ itself. It is efficient in the sense that it attains the Cramér–Rao bound

$$\sigma^2 d = \mathbb{E}_\mu[\|X - \mu\|_d^2] = \inf_Z \mathbb{E}_\mu[\|Z - \mu\|_d^2], \qquad \mu \in \mathbb{R}^d,$$

over all unbiased mean estimators $Z$ satisfying $\mathbb{E}_\mu[Z] = \mu$, for all $\mu \in \mathbb{R}^d$, where $\|\cdot\|_d$ denotes the Euclidean norm on $\mathbb{R}^d$.

In [5], James and Stein have constructed superefficient estimators for the mean of $X \in \mathbb{R}^d$, of the form

$$\left(1 - \frac{d-2}{\|X\|_d^2}\right) X,$$

whose risk is lower than the Cramér–Rao bound $\sigma^2 d$ in dimension $d \geq 3$.

The problem of drift estimation for Gaussian processes is of interest in several fields of application. For example, the process $(X_t)_{t \in [0,T]}$ can be interpreted as an observed output signal with decomposition

$$X_t = u_t + X_t^u, \qquad t \in [0,T],$$

where the drift $(u_t)_{t\in[0,T]}$ is viewed as an input signal to be estimated and perturbed by a centered Gaussian noise $(X^u_t)_{t\in[0,T]}$; see, for example, [4], Chapter VII. Such results find applications in, for example, telecommunication (additive Gaussian channels) and finance (identification of market trends).

Berger and Wolpert [2, 11], have constructed estimators of a James–Stein type for the drift of a Gaussian process $(X_t)_{t\in[0,T]}$ by applying the James–Stein procedure to the independent Gaussian random variables appearing in the Karhunen–Loève expansion of the process. In this context, $\hat{u} := (X_t)_{t\in\mathbb{R}_+}$ is seen as a minimax estimator of its own drift $(u_t)_{t\in\mathbb{R}_+}$.

Stein [10] has shown that the James–Stein estimators on $\mathbb{R}^d$ could be extended to a wider family of estimators, using integration by parts for Gaussian measures. Let us briefly recall Stein's argument, which relies on integration by parts with respect to the Gaussian density and on the properties of superharmonic functionals for the Laplacian on $\mathbb{R}^d$. Given an estimator of $\mu \in \mathbb{R}^d$ of the form $X + g(X)$, where $g:\mathbb{R}^d \to \mathbb{R}^d$ is sufficiently smooth, and applying the integration by parts formula

$$(1.1) \qquad \mathbb{E}_\mu[(X_i - \mu_i)g_i(X)] = \sigma^2 \mathbb{E}_\mu[\partial_i g_i(X)]$$

to $g = \sigma^2 \nabla \log f = \sigma^2(\partial_1 \log f, \ldots, \partial_d \log f)$, one obtains

$$\mathbb{E}_\mu[\|X + \sigma^2 \nabla \log f(X) - \mu\|_d^2] = \sigma^2 d + 4\sigma^4 \sum_{i=1}^d \mathbb{E}_\mu\left[\frac{\partial_i^2 \sqrt{f}(X)}{\sqrt{f}(X)}\right],$$

that is, $X + \sigma^2 \nabla \log f(X)$ is a superefficient estimator if

$$\sum_{i=1}^d \partial_i^2 \sqrt{f}(x) < 0, \qquad dx\text{-a.e.},$$

which is possible if $d \geq 3$. In this case, $X + \sigma^2 \nabla \log f(X)$ improves in the mean square sense over the efficient estimator $\hat{u}$ which attains the Cramér–Rao bound $\sigma^2 d$ on unbiased estimators of $\mu$.

In this paper we present an extension of Stein's argument to an infinite-dimensional setting using the Malliavin integration by parts formula, with application to the construction of Stein type estimators for the drift of a Gaussian process $(X_t)_{t\in[0,T]}$. Our approach applies to Gaussian processes such as Volterra processes and fractional Brownian motions. It also extends the results of Berger and Wolpert [2] in the same way that the construction of Stein [10] extends that of James and Stein [5], and this allows us to recover the estimators of the James–Stein type introduced by Berger and Wolpert [2] as particular cases. Here we replace the Stein equation (1.1) with the integration by parts formula of the Malliavin calculus on Gaussian space. Our estimators are given by processes of the form

$$X_t + D_t \log F, \qquad t \in [0,T],$$

where $F$ is a positive superharmonic random variable on Gaussian space and $D_t$ is the Malliavin derivative indexed by $t \in [0, T]$. In contrast to the minimax estimator $\hat{u}$, such estimators are not only biased but also anticipating with respect to the Brownian filtration $(\mathcal{F}_t)_{t \in [0,T]}$. This, however, poses no problem when one has access to complete paths from time $0$ to $T$.

When the variance of $X$ is large it can be shown that the percentage gain of this estimator is at least equal to the universal constant

$$(1.2) \qquad \frac{16}{\pi^4} \int_{\mathbb{R}^4} e^{-(x^2+y^2+z^2+r^2)/2} \frac{dx\,dy\,dz\,dr}{x^2 + 9y^2 + 25z^2 + 49r^2},$$

which approximately represents $11.38\%$; see (4.4) below.

We proceed as follows. Section 2 deals with notation and preliminaries on the representation of Gaussian processes. In Section 3 we show the estimator $\hat{u} := (X_t)_{t \in [0,T]}$ is minimax and in the case of independent increments we prove a Cramér–Rao bound over all unbiased drift estimators, which is attained by $\hat{u}$. Then we turn to the construction of superefficient drift estimators using the Malliavin gradient operator which show, as in the finite dimensional case, that the minimax estimator $\hat{u}$ is not admissible. We also construct nonnegative superharmonic functionals using cylindrical functionals on Gaussian space and show that the James–Stein estimators of Berger and Wolpert [2] can be recovered as particular cases of our approach. Examples and numerical simulations for the gain of such estimators are presented in Section 4 for the estimation of a deterministic drift. Proofs of the main results are provided in Section 5 using the Girsanov theorem, functional Bayes estimators and the Malliavin calculus.

**2. Notation.** Let $T > 0$. Consider a real-valued centered Gaussian process $(X_t)_{t \in [0,T]}$ with covariance function

$$\gamma(s,t) = \mathbb{E}[X_s X_t], \qquad s, t \in [0, T],$$

on a probability space $(\Omega, \mathcal{F}, \mathbb{P})$, where $\mathcal{F}$ is the $\sigma$-algebra generated by $X$.

We choose to represent $(X_t)_{t \in [0,T]}$ as an isonormal Gaussian process on the real separable Hilbert space $H$ generated by the functions $\chi_t(s) = \min(s,t)$, $s, t \in [0, T]$, with the scalar product $\langle \cdot, \cdot \rangle_H$ and norm $\|\cdot\|_H$ defined by

$$\langle \chi_t, \chi_s \rangle_H := \gamma(s,t), \qquad s, t \in [0, T].$$

In this case $X$ is viewed as an isometry $X : H \to L^2(\Omega, \mathcal{F}, P)$ with

$$X(\chi_t) := X_t, \qquad t \in [0, T],$$

and $\{X(h) : h \in H\}$ is a family of centered Gaussian random variables satisfying

$$\mathbb{E}[X(h)X(g)] = \langle h, g \rangle_H, \qquad h, g \in H.$$

We will assume in addition that $\gamma(s,t)$ has the form

$$\gamma(s,t) = \int_0^{s \wedge t} K(t,r)K(s,r)\,dr, \qquad s,t \in [0,T],$$

where $K(\cdot,\cdot)$ is a deterministic kernel and

$$(Kh)(t) := \int_0^t K(t,s)\dot{h}(s)\,ds$$

is differentiable in $t \in [0,T]$. In this case the scalar product in $H$ satisfies

$$\langle h, g \rangle_H = \langle K^*h, K^*g \rangle = \langle h, \Gamma g \rangle,$$

where $\Gamma = KK^*$ and $K^*$ is the adjoint of $K$ with respect to

$$\langle h, g \rangle := \langle \dot{h}, \dot{g} \rangle_{L^2([0,T],dt)};$$

see [1]. Moreover, for any orthonormal basis $(h_k)_{k \in \mathbb{N}}$ of $H$, we have the expansion

$$
\begin{aligned}
X_t &= \sum_{k=0}^{\infty} \langle \chi_t, h_k \rangle_H X(h_k) \\
&= \sum_{k=0}^{\infty} \langle 1_{[0,t]}, \dot{\Gamma} h_k \rangle_{L^2([0,T],dt)} X(h_k) \\
&= \sum_{k=0}^{\infty} \Gamma h_k(t) X(h_k),
\end{aligned}
\tag{2.1}
$$

$t \in [0,T]$, and the representation

$$X_t = \int_0^t K(t,s)\,dW_s, \qquad t \in [0,T],$$

where $(W_s)_{s \in [0,T]}$ is a standard Brownian motion; see [1].

Note that there exists other constructions of Gaussian processes leading to different series decompositions for $X_t$, for example, Berger and Wolpert [2, 11] use Karhunen–Loève expansions; see [8] for a formulation of our results in that setting.

## 3. Main results.

*Efficient drift estimator.* Recall that the classical linear parametric estimation problem for the drift of a diffusion consists in estimating the coefficient $\theta$ appearing in

$$d\xi_t = \theta a_t(\xi_t)\,dt + dY_t, \qquad \xi_0 = 0,$$

with a maximum likelihood estimator $\hat{\theta}_T$ given by

$$\hat{\theta}_T = \frac{\int_0^T a_t(\xi_t)\,d\xi_t}{\int_0^T a_t^2(\xi_t)\,dt}; \tag{3.1}$$

see [6, 9] in case $(Y_t)_{t\in\mathbb{R}_+}$ is a Brownian motion.

In this paper we consider the nonparametric functional estimation of the drift of a one-dimensional Gaussian process $(X_t)_{t\in\mathbb{R}_+}$ with decomposition

$$dX_t = \dot{u}_t\,dt + dX_t^u, \tag{3.2}$$

where $(u_t)_{t\in[0,T]}$ is an adapted process of the form

$$u_t = \int_0^t \dot{u}_s\,ds, \qquad t\in[0,T],\ \dot{u}\in L^2(\Omega\times[0,T]),$$

and $(X_t^u)_{t\in\mathbb{R}_+}$ is a centered Gaussian process under a probability $\mathbb{P}_u$ which is the translation of $\mathbb{P}$ on $\Omega$ by $u$. The expectation under $\mathbb{P}_u$ will be denoted by $\mathbb{E}_u$.

DEFINITION 3.1.  A drift estimator $\xi$ is called unbiased if

$$\mathbb{E}_u[\xi_t] = \mathbb{E}_u[u_t], \qquad t\in[0,T],$$

for all square-integrable adapted process $(u_t)_{t\in[0,T]}$. It is called adapted if the process $(\xi_t)_{t\in[0,T]}$ is adapted to the filtration $(\mathcal{F}_t)_{t\in[0,T]}$ generated by $(X_t^u)_{t\in[0,T]}$.

Here, the canonical process $(X_t)_{t\in[0,T]}$ will be considered as an unbiased estimator of its drift $(u_t)_{t\in[0,T]}$ under $\mathbb{P}_u$, with risk defined as

$$\mathrm{R}(\gamma,\mu,\hat{u}) := \mathbb{E}_u[\|X-u\|_{L^2([0,T],d\mu)}^2]$$

$$= \int_0^T \mathbb{E}_u[|X_t^u|^2]\mu(dt)$$

$$= \int_0^T \gamma(t,t)\mu(dt),$$

where $\mu$ is a finite weighting Borel measure on $[0,T]$. In case $u$ is constrained to have the form $u_t = \theta t$, $t\in[0,T]$, $\theta\in\mathbb{R}$, our estimator $\hat{u}$ satisfies $\hat{u}_T = \hat{\theta}_T T$, where $\hat{\theta}_T$ is given by (3.1), $T>0$, with the asymptotics $\hat{\theta}_T \to \theta$ in probability as $T$ tends to infinity.

In the next proposition we note that, as in the finite dimensional Gaussian case, the estimator $\hat{u} = (X_t)_{t\in[0,T]}$ is minimax. Note that no adaptedness condition is imposed on $\xi$ in the infimum (3.3).

PROPOSITION 3.2. *The estimator $\hat{u} = (X_t)_{t \in [0,T]}$ is minimax. For all $u \in \Omega$, we have*

$$
\begin{aligned}
\mathrm{R}(\gamma, \mu, \hat{u}) &= \mathbb{E}_u \left[ \int_0^T |X_t - u_t|^2 \mu(dt) \right] \\
&= \inf_\xi \sup_{v \in \Omega} \mathbb{E}_v \left[ \int_0^T |\xi_t - v_t|^2 \mu(dt) \right].
\end{aligned}
\tag{3.3}
$$

This proposition is proved in Section 5 using functional Bayes estimators.

We close this section with a more precise statement concerning the efficiency of the estimator $\hat{u} = (X_t)_{t \in [0,T]}$ in the case where $(X_t)_{t \in [0,T]}$ has independent increments, that is, $\gamma(s,t)$ is of the form

$$\gamma(s,t) = \int_0^{s \wedge t} \sigma_u^2 \, du,$$

where $\sigma \in L^2([0,T], dt)$ is an a.e. nonvanishing function with

$$
\mathbb{E} \left[ \int_0^T \frac{\dot{u}_s^2}{\sigma_s^2} ds \right] < \infty.
\tag{3.4}
$$

In this case the following proposition allows us to compute a Cramér–Rao bound attained by $\hat{u}$, and shows that $\hat{u} = (X_t)_{t \in [0,T]}$ is an efficient estimator.

PROPOSITION 3.3. *Cramér–Rao inequality. For any unbiased and adapted estimator $\xi$ of $u$ we have*

$$
\mathbb{E}_u \left[ \int_0^T |\xi_t - u_t|^2 \mu(dt) \right] \geq \mathrm{R}(\sigma, \mu, \hat{u}),
\tag{3.5}
$$

*where $u \in L^2(\Omega \times [0,T], \mathbb{P}_u \otimes \mu)$ is adapted and the Cramér–Rao type bound*

$$\mathrm{R}(\sigma, \mu, \hat{u}) := \int_0^T \int_0^t \sigma_s^2 \, ds \, \mu(dt)$$

*is independent of $u$ and is attained by the efficient estimator $\hat{u} = X$.*

This proposition is proved in Section 5 using the Girsanov theorem.

*Superefficient drift estimators.* Next we construct a family of superefficient estimators of $u$ of the form $X + \xi$, whose mean square error is strictly smaller than the minimax risk $\mathrm{R}(\gamma, \mu, \hat{u})$ of Proposition 3.2 when $\xi \in L^2([0,T] \times \Omega, \mathbb{P}_u \otimes \mu)$ is a suitably chosen stochastic process. This estimator will be biased and anticipating with respect to the Brownian filtration. In the next theorem we follow Stein's argument which uses integration by parts, but we replace (1.1) by the duality relation (5.5) between the gradient and divergence operators on Gaussian space.

Before turning to the main result, we need to introduce some elements of analysis and Malliavin calculus on Gaussian space, see, for example, [7]. Fix $(h_n)_{n\geq 1}$ a total subset of $H$, and let $\mathcal{S}$ denote the space of cylindrical functionals of the form

$$(3.6) \qquad F = f_n(X^u(h_1), \ldots, X^u(h_n)),$$

where $f_n$ is in the space of infinitely differentiable rapidly decreasing functions on $\mathbb{R}^n$, $n \geq 1$.

DEFINITION 3.4. The Malliavin derivative $D$ is defined as

$$D_t F = \sum_{i=1}^n \Gamma h_i(t) \partial_i f_n(X^u(h_1), \ldots, X^u(h_n)), \qquad t \in [0, T],$$

for $F \in \mathcal{S}$ of the form (3.6).

It is known that $D$ is closable (cf. Proposition 1.2.1 of [7]) and its closed domain will be denoted by $\mathrm{Dom}(D)$.

DEFINITION 3.5. We define the Laplacian $\Delta$ by

$$\Delta F = \mathrm{trace}_{L^2([0,T],d\mu)^{\otimes 2}} DDF = \int_0^T D_t D_t F \mu(dt)$$

on the space $\mathrm{Dom}(\Delta)$ made of all $F \in \mathrm{Dom}(D)$ such that $D_t F \in \mathrm{Dom}(D)$, $t \in [0, T]$, and $(D_t D_t F)_{t \in [0,T]} \in L^2([0, T], d\mu)$, $\mathbb{P}$-a.s.

For $F \in \mathcal{S}$ of the form (3.6), we have

$$\Delta F = \sum_{i,j=1}^n \langle \Gamma h_i, \Gamma h_j \rangle_{L^2([0,T],d\mu)} \partial_i \partial_j f_n(X^u(h_1), \ldots, X^u(h_n)),$$

and it can be easily shown that the operator $\Delta$ is closable; see, for example, [8]. We will say that a random variable $F$ in $\mathrm{Dom}(\Delta)$ is $\Delta$-superharmonic on $\Omega$ if

$$(3.7) \qquad \Delta F(\omega) \leq 0, \qquad \mathbb{P}(d\omega)\text{-a.s.}$$

The next theorem is our main result on the construction of superefficient estimators. It is proved in Section 5 using the Malliavin calculus on Gaussian space.

THEOREM 3.6. (i) *Unbiased risk estimate.* For any $\xi \in L^2(\Omega \times [0,T], \mathbb{P}_u \otimes \mu)$ such that $\xi_t \in \mathrm{Dom}(D)$, $t \in [0,T]$, and $(D_t \xi_t)_{t \in [0,T]} \in L^1(\Omega \times [0,T], \mathbb{P}_u \otimes \mu)$, we have

$$
\mathbb{E}_u[\|X + \xi - u\|^2_{L^2([0,T], d\mu)}] \tag{3.8}
$$
$$
= \mathrm{R}(\gamma, \mu, \hat{u}) + \|\xi\|^2_{L^2(\Omega \times [0,T], \mathbb{P}_u \otimes \mu)} + 2\mathbb{E}_u\left[\int_0^T D_t \xi_t \mu(dt)\right].
$$

(ii) *Stein-type estimator.* For any $\mathbb{P}$-*a.s.* positive random variable $F \in \mathrm{Dom}(D)$ such that $D_t F \in \mathrm{Dom}(D)$, $t \in [0,T]$, and $(D_t D_t F)_{t \in [0,T]} \in L^1(\Omega \times [0,T], \mathbb{P}_u \otimes \mu)$, we have

$$
\mathbb{E}_u[\|X + D\log F - u\|^2_{L^2([0,T], d\mu)}]
$$
$$
(3.9) \qquad = \mathrm{R}(\gamma, \mu, \hat{u}) - \mathbb{E}_u[\|D\log F\|^2_{L^2([0,T], d\mu)}] + 2\mathbb{E}_u\left[\frac{\Delta F}{F}\right]
$$
$$
= \mathrm{R}(\gamma, \mu, \hat{u}) + 4\mathbb{E}_u\left[\frac{\Delta \sqrt{F}}{\sqrt{F}}\right].
$$

REMARKS.

(a) The $\Delta$-superharmonicity of $F$ is sufficient but not necessary for $X + D\log F$ to be superefficient, namely, it can be replaced by the $\Delta$-superharmonicity of $\sqrt{F}$, which is a weaker assumption; see [3] in the finite dimensional case.

(b) As in [10], the superefficient estimators constructed in this way are minimax in the sense that, from Proposition 3.2 and Theorem 3.6(ii), for all $u \in H$ we have,

$$
\mathbb{E}_u[\|X + D\log F - u\|^2_{L^2([0,T], d\mu)}] < \mathrm{R}(\gamma, \mu, \hat{u}) = \inf_{\xi} \sup_{v \in \Omega} \mathbb{E}_v\left[\int_0^T |\xi_t - v_t|^2 dt\right],
$$

provided $\Delta \sqrt{F} < 0$ on a set of strictly positive $\mathbb{P}$-measure, thus showing that the minimax estimator $\hat{u} = (X_t)_{t \in [0,T]}$ is inadmissible.

We now turn to some examples of nonnegative superharmonic functionals with respect to the Laplacian $\Delta$, which satisfy the hypotheses of Theorem 3.6 and will be used in the numerical applications of Section 4. We assume that $(\Gamma h_k)_{k \geq 1}$ is orthogonal in $L^2([0,T], d\mu)$, and we let

$$
\lambda_k = \|\Gamma h_k\|_{L^2([0,T], d\mu)}, \qquad k \geq 1.
$$

For example, in case $\mu(dt) = dt$ the sequence $(h_k)_{k \geq 1}$ can be realized as the solution of the eigenvalue problem

$$
(3.10) \qquad \Gamma h_k = -\lambda_k^2 \ddot{h}_k, \qquad \dot{h}_k(T) = 0, \qquad k \geq 1.
$$

Superharmonic functionals on Gaussian space can be constructed as cylindrical functionals, by composition with finite-dimensional functions. From the expression

$$\Delta F = \sum_{i=1}^n \partial_i^2 f_n(\lambda_1^{-1} X^u(h_1), \ldots, \lambda_n^{-1} X^u(h_n)),$$

of $\Delta$ on cylindrical functionals, we check that

$$F = f_n(\lambda_1^{-1} X^u(h_1), \ldots, \lambda_n^{-1} X^u(h_n))$$

is superharmonic on $\Omega$ if and only if $f_n$ is superharmonic on $\mathbb{R}^n$. Note that when $u \in H$ is deterministic, any superharmonic functional of the form

$$f_n(\lambda_1^{-1} X^u(h_1), \ldots, \lambda_n^{-1} X^u(h_n))$$

can be replaced with

$$f_n(\lambda_1^{-1} X(h_1), \ldots, \lambda_n^{-1} X(h_n)),$$

which retains the same harmonicity property since $u$ is deterministic, and can be directly computed from an observation of $X$.

Given $a \in \mathbb{R}$ and $b \in \mathbb{R}^n$, let $f_{n,a,b} : \mathbb{R}^n \to \mathbb{R}$ be defined as

$$f_{n,a,b}(x_1, \ldots, x_n) = \|x + b\|^a = ((x_1 + b_1)^2 + \cdots + (x_n + b_n)^2)^{a/2},$$

then $\sqrt{f_{n,a,b}}$ is superharmonic on $\mathbb{R}^n$, $n \geq 3$, if and only if $a \in [4 - 2n, 0]$. Letting

$$F_{n,a,b} = f_{n,a,b}(\lambda_1^{-1} X^u(h_1), \ldots, \lambda_n^{-1} X^u(h_n)),$$

we have

$$D_t \log F_{n,a,b} = a \sum_{i=1}^n \frac{\lambda_i^{-1} \Gamma h_i(t)(b_i + \lambda_i^{-1} X^u(h_i))}{|b_1 + \lambda_1^{-1} X^u(h_1)|^2 + \cdots + |b_n + \lambda_n^{-1} X^u(h_n)|^2}$$

and

$$\Delta \sqrt{F_{n,a,b}} = \sum_{i=1}^n \partial_i^2 \sqrt{f_{n,a,b}}(\lambda_1^{-1} X^u(h_1), \ldots, \lambda_n^{-1} X^u(h_n)),$$

since $(\Gamma h_k)_{k \geq 1}$ is orthogonal in $L^2([0, T], d\mu)$. Hence,

$$(3.11) \quad \frac{\Delta \sqrt{F_{n,a,b}}}{\sqrt{F_{n,a,b}}} = \frac{a(n - 2 + a/2)/2}{|b_1 + \lambda_1^{-1} X^u(h_1)|^2 + \cdots + |b_n + \lambda_n^{-1} X^u(h_n)|^2}$$

is negative if $4 - 2n \leq a \leq 0$. On the other hand,

$$(3.12) \quad \frac{\Delta F_{n,a,b}}{F_{n,a,b}} = \frac{a(n + a - 2)}{|b_1 + \lambda_1^{-1} X^u(h_1)|^2 + \cdots + |b_n + \lambda_n^{-1} X^u(h_n)|^2}$$

is negative for $a \in (2-n, 0]$ and vanishes for $a = 2-n$. Taking $b_i = \lambda_i^{-1} \langle u, h_i \rangle$, $i = 1, \ldots, n$, we have

$$D_t \log F_{n,2-n,b} = -(n-2) \frac{[\Pi_n X]_t}{\|\Pi_n X\|^2_{L^2([0,T],dt)}},$$

where $\Pi_n$ denotes the orthogonal projection

$$\Pi_n X(t) := \sum_{k=1}^{n} \lambda_k^{-1} X(h_k) \Gamma h_k(t) = \sum_{k=1}^{n} \lambda_k^{-1} (b_k + \lambda_k^{-1} X^u(h_k)) \Gamma h_k(t).$$

The resulting estimator

$$X_t + D_t \log F_{n,2-n,b} = X_t - (n-2) \frac{[\Pi_n X]_t}{\|\Pi_n X\|^2_{L^2([0,T],dt)}}, \qquad t \in [0,T],$$

is of a James–Stein type, but it is not a shrinkage operator.

From (3.12) we have $\Delta F_{n,2-n,b} = 0$, hence (3.9) and (3.11) show that

$$\|D \log F_{n,2-n,b}\|^2_{L^2([0,T] \times \Omega, \mathbb{P}_u \otimes dt)}$$

$$= -4 \mathbb{E}_u \left[ \frac{\Delta \sqrt{F_{n,2-n,b}}}{\sqrt{F_{n,2-n,b}}} \right]$$

(3.13)

$$= (n-2)^2 \mathbb{E}_u \left[ \frac{1}{|\lambda_1^{-1} X(h_1)|^2 + \cdots + |\lambda_n^{-1} X(h_n)|^2} \right]$$

$$= (n-2)^2 \mathbb{E}_u[\|\Pi_n X\|^{-2}_{L^2([0,T],dt)}],$$

and the risk of $X + D \log F_{n,2-n,b}$ is

$$\mathbb{E}_u[\|X + D \log F_{n,2-n,b} - u\|^2_{L^2([0,T],dt)}]$$

$$= \mathrm{R}(\gamma, \mu, \hat{u}) - (n-2)^2 \mathbb{E}_u[\|\Pi_n X\|^{-2}_{L^2([0,T],dt)}].$$

**4. Numerical application.** In this section we present numerical simulations which allow us to measure the efficiency of our estimators. We use the superharmonic functionals constructed as cylindrical functionals in the previous section, and we assume that $u \in H$ is deterministic.

We work in the case of independent increments and we additionally assume that $\sigma_t = \sigma$ is constant, $t \in [0,T]$, that is, $(X_t)_{t \in [0,T]}$ is a Brownian motion with variance $\sigma^2$, $\Gamma h(t) = \sigma^2 h(t)$, $t \in [0,T]$, and

$$\mathrm{R}(\sigma, \mu, \hat{u}) = \frac{\sigma^2}{2} T^2.$$

Letting

$$h_n(t) = \frac{\sqrt{2T}}{\sigma \pi (n - 1/2)} \sin\left( \left(n - \frac{1}{2}\right) \frac{\pi t}{T} \right), \qquad t \in [0,T], \ n \geq 1,$$

that is,

$$\dot{h}_n(t) = \frac{1}{\sigma}\sqrt{\frac{2}{T}}\cos\left(\left(n-\frac{1}{2}\right)\frac{\pi t}{T}\right), \qquad t \in [0,T],\ n \geq 1,$$

provides an orthonormal basis $(h_n)_{n\geq 1}$ of $H$ such that $(\Gamma h_k)_{k\geq 1}$ is orthogonal in $L^2([0,T], dt)$, with

$$\lambda_n = \frac{\sigma T}{\pi(n-1/2)}, \qquad n \geq 1,$$

and satisfies (3.10). The Stein estimator $X + D\log F$ of $u$ will be given by

$$D_t \log F_{n,2-n,b} = -(n-2)\sqrt{\frac{2}{T}}\sum_{k=1}^{n}\frac{X(h_k)}{|\lambda_1^{-1}X(h_1)|^2 + \cdots + |\lambda_n^{-1}X(h_n)|^2}$$
$$\times \sin\left(\left(k-\frac{1}{2}\right)\frac{\pi t}{T}\right).$$

For simulation purposes we will represent the (nondrifted) Brownian motion $(X_t^u)_{t\in[0,T]}$ via the Paley–Wiener expansion

$$(4.1) \qquad X_t^u = \sigma^2 \sum_{n=1}^{\infty} \eta_n h_n(t) = \sigma\frac{\sqrt{2T}}{\pi}\sum_{n=1}^{\infty} \eta_n \frac{\sin((n-1/2)(\pi t)/T)}{(n-1/2)},$$

where $(\eta_n)_{n\geq 1}$ are independent standard Gaussian random variables with unit variance under $\mathbb{P}_u$ and

$$\eta_n = \int_0^T \dot{h}_n(s)\,dX_s^u, \qquad n \geq 1.$$

In this case we have

$$(4.2) \qquad D_t \log F_{n,2-n,b} = -(n-2)\sqrt{\frac{2}{T}}\sum_{k=1}^{n}\frac{\eta_k + \langle u, h_k\rangle}{\sum_{l=1}^{n} \lambda_l^{-2}(\eta_l + \langle u, h_l\rangle)^2}$$
$$\times \sin\left(\left(k-\frac{1}{2}\right)\frac{\pi t}{T}\right).$$

The gain of the superefficient estimator $X + D\log F_{n,2-n,b}$ compared to the minimax and efficient estimator $\hat{u}$ is given by

$$G(u,\sigma,T,n) := -\frac{4}{\mathrm{R}(\sigma,\mu,\hat{u})}\mathbb{E}_u\left[\frac{\Delta\sqrt{F_{n,2-n,b}}}{\sqrt{F_{n,2-n,b}}}\right], \qquad n \geq 3.$$

From (3.9) and (3.13) we have

$$(4.3) \qquad G(u,\sigma,T,n) = 2(n-2)^2 \mathbb{E}\left[\left(\sum_{l=1}^{n}\left(\pi\left(l-\frac{1}{2}\right)(\eta_l + \langle u, h_l\rangle)\right)^2\right)^{-1}\right],$$

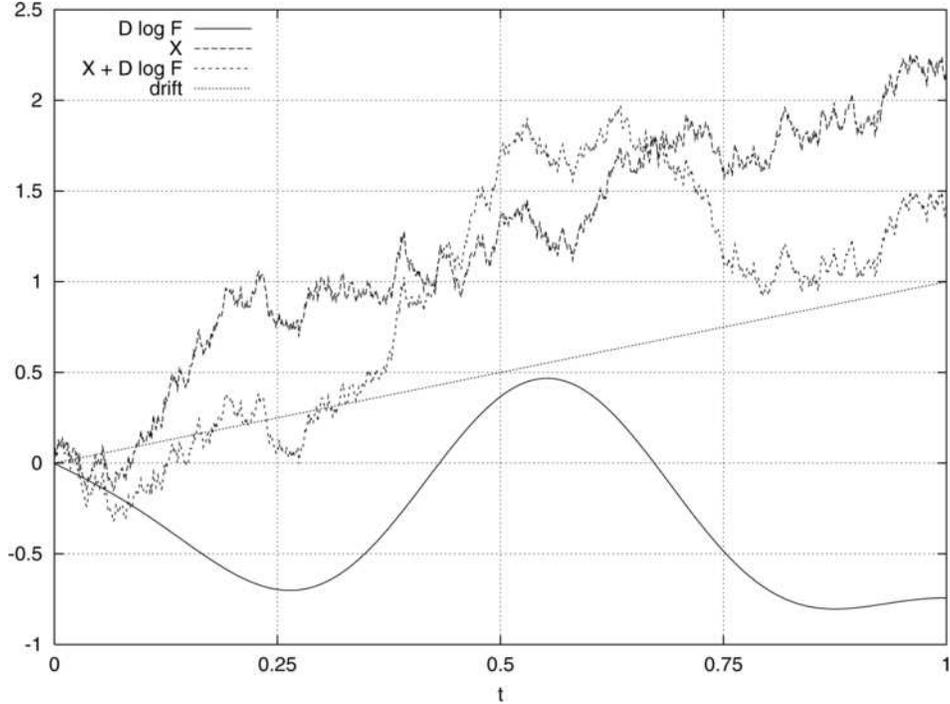

Fig. 1.  $u(t) = t, t \in [0, T]; \ n = 5.$

hence,

$$(4.4) \qquad \lim_{\sigma \to \infty} G(u, \sigma, T, n) = (n-2)^2 \frac{8}{\pi^2} \mathbb{E}\left[\left(\sum_{l=1}^{n} (2l-1)^2 \eta_l^2\right)^{-1}\right].$$

The quantity (4.4) can be evaluated as a Gaussian integral to yield (1.2). Unlike in the classical Stein method, we stress that here $n$ becomes a free parameter and there is some interest in determining the values of $n$ which yield the best performance. The next proposition is proved by standard arguments (see [8] for details), and is illustrated by Figure 2.

PROPOSITION 4.1.  *For all $\sigma, T > 0$, and $u \in H$, we have*

$$G(u, \sigma, T, n) \simeq \frac{6}{n\pi^2}$$

*as $n$ goes to infinity.*

In the sequel we choose $u_t = \alpha t$, $t \in [0, T]$, $\alpha \in \mathbb{R}$. Figure 1 gives a sample path representation of the process $X + D \log F$.

In this case, from (4.3) we have

$$G(\alpha,\sigma,T,n) = 2(n-2)^2 \mathbb{E}\left[\left(\sum_{l=1}^{n}\left(\pi\left(l-\frac{1}{2}\right)\eta_l - \alpha\frac{\sqrt{2T}}{\sigma}(-1)^l\right)^2\right)^{-1}\right],$$

from which it follows that $G(\alpha,\sigma,T,n)$ converges to

$$(n-2)^2\frac{8}{\pi^2}\mathbb{E}\left[\left(\sum_{l=1}^{n}(2l-1)^2\eta_l^2\right)^{-1}\right],$$

when $\alpha^{-2}\sigma^2/T$ tends to infinity, and is equivalent to

$$\left(1-\frac{2}{n}\right)^2\frac{\sigma^2}{\alpha^2 T}$$

as $\alpha^{-2}\sigma^2/T$ tends to 0. Figure 2 represents the gain in percentage of the superefficient estimator $X + \sigma^2 D \log F_{n,2-n,b}$ compared to the efficient estimator $\hat{u}$ using Monte Carlo simulations, that is, we represent $100 \times G(\alpha,\sigma,T,n)$ as a function of $n \geq 3$. An optimal value

$$n_{\text{opt}} = \arg\max\{G(\alpha,\sigma,T,n) : n \geq 3\}$$

of $n$ exists in general and is equal to 4 when $\alpha = \sigma = T = 1$. Figure 3 shows the variation of the gain as a function of $n$ and $T$ for $\alpha = \sigma = 1$. Figure 4 represents the variation of the gain as a function of $n$ and $\sigma$.

## 5. Proofs of the main results.

*Minimaxity.* In order to prove the minimaxity of the estimator $\hat{u} = X$, we will need the notion of Bayes estimator. We will make use of the next lemma which is classical in the framework of Gaussian filtering. We say that a Gaussian process $Z$ has covariance $\Gamma$ and drift $v \in \Omega$ if $Z - v$ is a centered Gaussian process with covariance $\Gamma$.

LEMMA 5.1. *Let $Z$ be a Gaussian process with covariance operator $\Gamma_\tau$ and drift $v \in H$, and assume that $X$ is a Gaussian process with drift $Z$ and covariance operator $\Gamma$ given $Z$. Then, conditionally to $X$, $Z$ has drift*

$$f \mapsto \langle f, (\Gamma + \Gamma_\tau)^{-1}\Gamma v\rangle + X((\Gamma+\Gamma_\tau)^{-1}\Gamma_\tau f)$$

*and covariance*

$$\Gamma_\tau(\Gamma + \Gamma_\tau)^{-1}\Gamma.$$

Note that, unlike in Proposition 3.3, no adaptedness or unbiasedness restriction is made on $\xi$ in the infimum taken in (5.2) below.

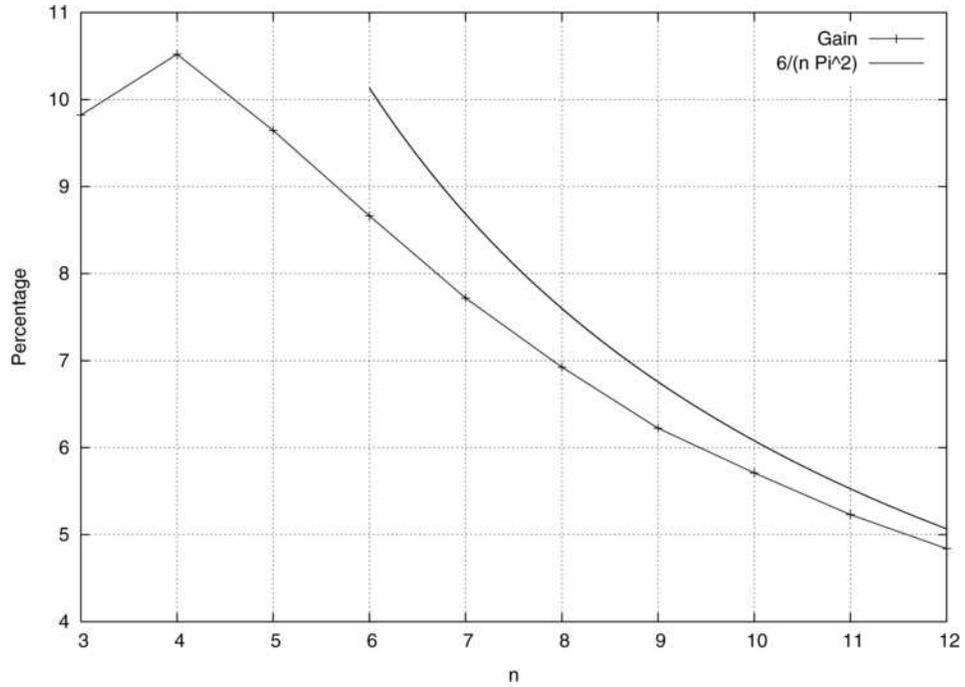

Fig. 2. *Percentage gain as a function of n for 10000 samples and $\alpha = \sigma = T = 1$.*

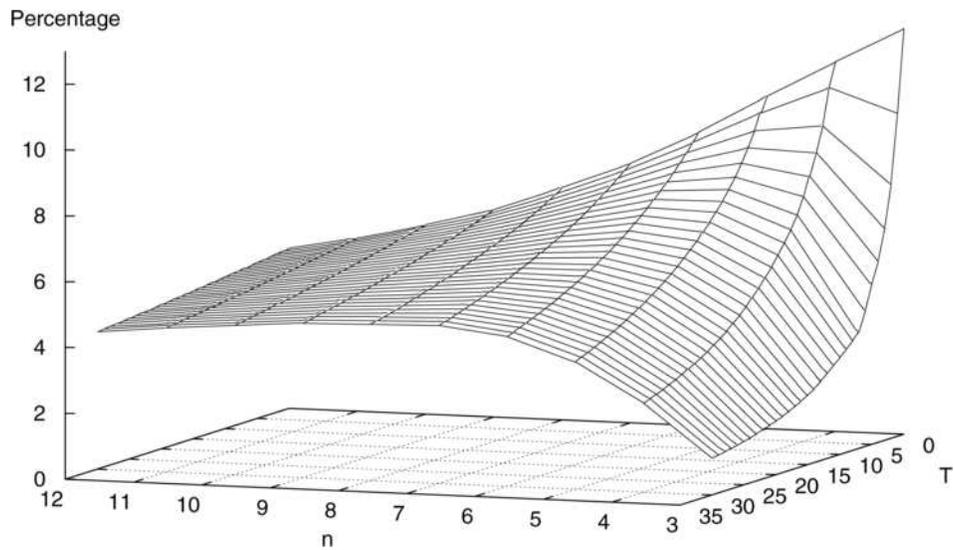

Fig. 3. *Gain as a function of n and T.*

PROPOSITION 5.2. *Bayes estimator.* Let $\mathbb{P}_v^\tau$ denote the Gaussian distribution on $\Omega$ with covariance operator $\Gamma_\tau$ and drift $v \in H$. The Bayes risk

$$\text{(5.1)} \qquad \int_\Omega \mathbb{E}_z \left[ \int_0^T |\xi_t - z_t|^2 \mu(dt) \right] d\mathbb{P}_v^\tau(z)$$

*of any drift estimator* $(\xi_t)_{t \in [0,T]}$ *on* $\Omega$ *under the prior distribution* $\mathbb{P}_v^\tau$ *is uniquely minimized by*

$$\xi_t^{\tau,v} := \langle \chi_t, (\Gamma_\tau + \Gamma)^{-1} \Gamma v \rangle + X((\Gamma_\tau + \Gamma)^{-1} \Gamma_\tau \chi_t), \qquad t \in [0, T],$$

*which has risk*

$$\text{(5.2)} \qquad \begin{aligned} &\int_0^T \langle \chi_t, \Gamma(\Gamma_\tau + \Gamma)^{-1} \Gamma_\tau \chi_t \rangle \mu(dt) \\ &= \inf_\xi \int_\Omega \mathbb{E}_z \left[ \int_0^T |\xi_t - z_t|^2 \mu(dt) \right] d\mathbb{P}_v^\tau(z). \end{aligned}$$

PROOF. Let $Z$ denote a Gaussian process with drift $v \in H$ and covariance $\Gamma_\tau$. From Lemma 5.1, $(Z_t)_{t \in [0,T]}$ has drift

$$t \mapsto \langle \chi_t, (\Gamma_\tau + \Gamma)^{-1} \Gamma v \rangle + X((\Gamma_\tau + \Gamma)^{-1} \Gamma_\tau \chi_t), \qquad t \in [0, T],$$

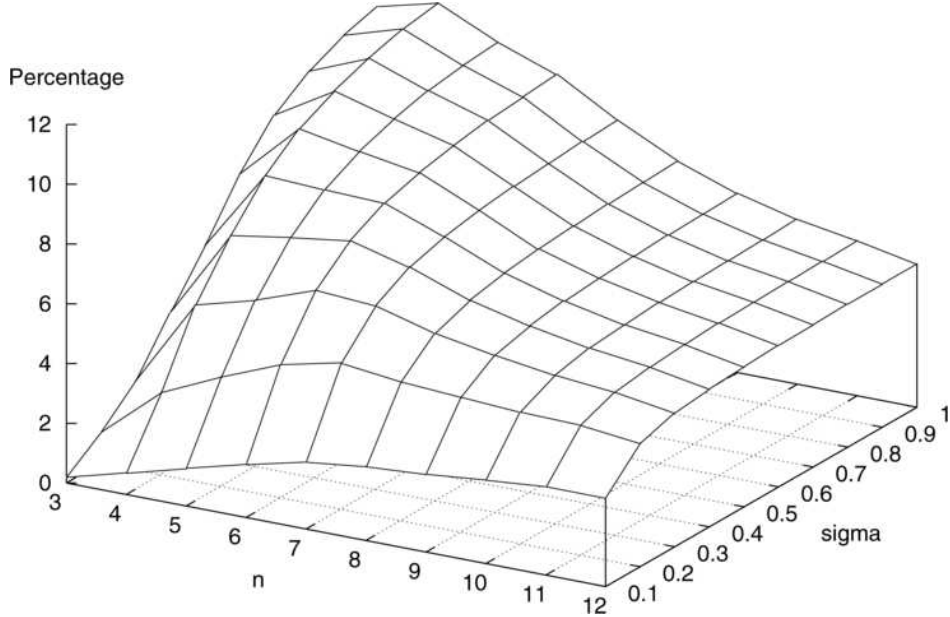

FIG. 4. *Gain as a function of $n$ and $\sigma$.*

and covariance $\Gamma_\tau(\Gamma_\tau + \Gamma)^{-1}\Gamma$ conditionally to $X$. Hence, the Bayes risk of an estimator $\xi$ under the prior distribution $\mathbb{P}_v^\tau$ is given by

$$\int_\Omega \mathbb{E}_z\left[\int_0^T |\xi_t - z_t|^2 \mu(dt)\right] d\mathbb{P}_v^\tau(z)$$

$$= \mathbb{E}\left[\mathbb{E}\left[\int_0^T |\xi_t - Z_t|^2 \mu(dt)\Big|X\right]\right]$$

$$= \mathbb{E}\left[\int_0^T |\xi_t - \mathbb{E}[Z_t \mid X]|^2 \mu(dt)\right] + \mathbb{E}\left[\int_0^T \mathrm{Var}(Z_t|X)\mu(dt)\right]$$

$$= \mathbb{E}\left[\int_0^T |\xi_t - \langle \chi_t, (\Gamma_\tau + \Gamma)^{-1}\Gamma v\rangle - X((\Gamma_\tau + \Gamma)^{-1}\Gamma_\tau \chi_t)|^2 \mu(dt)\right]$$

$$+ \int_0^T \langle \chi_t, \Gamma(\Gamma_\tau + \Gamma)^{-1}\Gamma_\tau \chi_t\rangle \mu(dt),$$

which is minimized by

$$\xi_t^{\tau,v} := \mathbb{E}[Z_t \mid X] = \langle \chi_t, (\Gamma_\tau + \Gamma)^{-1}\Gamma v\rangle - X((\Gamma_\tau + \Gamma)^{-1}\Gamma_\tau \chi_t),$$
$$t \in [0, T].$$

□

Clearly, $\xi^{\tau,v}$ is unique in the sense that it is the only estimator which minimizes the Bayes risk (5.1). This shows, in particular, that every $\xi^{\tau,v}$ is admissible in the sense that if an estimator $\xi$ satisfies

$$\mathbb{E}_z[\|\xi - z\|_{L^2([0,T],d\mu)}^2] \leq \mathbb{E}_z[\|\xi^{\tau,v} - z\|_{L^2([0,T],d\mu)}^2], \qquad z \in \Omega,$$

then

$$\int_\Omega \mathbb{E}_z[\|\xi - z\|_{L^2([0,T],d\mu)}^2] d\mathbb{P}_v^\tau(z) \leq \int_\Omega \mathbb{E}_z[\|\xi^{\tau,v} - z\|_{L^2([0,T],d\mu)}^2] d\mathbb{P}_v^\tau(z)$$

$$= \int_0^T \langle \chi_t, \Gamma(\Gamma_\tau + \Gamma)^{-1}\Gamma_\tau \chi_t\rangle \mu(dt),$$

hence,

(5.3) $$\int_\Omega \mathbb{E}_z[\|\xi - z\|_{L^2([0,T],d\mu)}^2] d\mathbb{P}_v^\tau(z) = \int_0^T \langle \chi_t, \Gamma(\Gamma_\tau + \Gamma)^{-1}\Gamma_\tau \chi_t\rangle \mu(dt),$$

and $\xi = \xi^{\tau,v}$ by Proposition 5.2.

PROOF OF PROPOSITION 3.2. We apply Proposition 5.2 in the case $\Gamma_\tau f(t) = \tau^2 f(t)$, $t \in [0, T]$. Clearly, taking $\xi = 0$ yields

$$\mathrm{R}(\gamma, \mu, \hat{u}) = \sup_{u \in \Omega} \mathbb{E}_u\left[\int_0^T |X_t - u_t|^2 \mu(dt)\right] \geq \inf_\xi \sup_{u \in \Omega} \mathbb{E}_u\left[\int_0^T |\xi_t - u_t|^2 \mu(dt)\right].$$

On the other hand, from Proposition 5.2, for all processes $\xi$ we have

$$\sup_{u \in \Omega} \mathbb{E}_u \left[ \int_0^T |\xi_t - u_t|^2 \mu(dt) \right] \geq \int_\Omega \mathbb{E}_z \left[ \int_0^T |\xi_t - z_t|^2 \mu(dt) \right] d\mathbb{P}_0^\tau(z)$$

$$\geq \int_0^T \langle \chi_t, (I + \Gamma/\tau^2)^{-1} \Gamma \chi_t \rangle \mu(dt),$$

for all $\tau > 0$, hence,

$$\inf_\xi \sup_{u \in H} \mathbb{E}_u \left[ \int_0^T |\xi_t - u_t|^2 \mu(dt) \right] \geq \int_0^T \langle \chi_t, \Gamma \chi_t \rangle \mu(dt) = \mathrm{R}(\gamma, \mu, \hat{u}). \qquad \square$$

*Cramér–Rao bound.*

PROOF OF THE CRAMÉR–RAO INEQUALITY PROPOSITION 3.3. By the Girsanov theorem, $\mathbb{P}_u$ is absolutely continuous with respect to $\mathbb{P}$ under the condition (3.4), with

$$d\mathbb{P}_u = \Lambda(u) \, d\mathbb{P},$$

where

$$\Lambda(u) := \exp\left( \int_0^T \frac{\dot{u}_s}{\sigma_s^2} dX_s - \frac{1}{2} \int_0^T \frac{\dot{u}_s^2}{\sigma_s^2} ds \right)$$

denotes the Girsanov–Cameron–Martin density, and the canonical process $(X_t)_{t \in [0,T]}$ becomes a continuous Gaussian semimartingale under $\mathbb{P}_u$, with quadratic variation $\sigma_t^2 \, dt$ and drift $\dot{u}_t \, dt$.

Since $\xi$ is unbiased, for all $\zeta \in H$, we have

$$\mathbb{E}_{u+\varepsilon\zeta}[\xi_t] = \mathbb{E}_{u+\varepsilon\zeta}[u_t + \varepsilon\zeta_t]$$
$$= \mathbb{E}_{u+\varepsilon\zeta}[u_t] + \varepsilon \mathbb{E}_{u+\varepsilon\zeta}[\zeta_t]$$
$$= \mathbb{E}_{u+\varepsilon\zeta}[u_t] + \varepsilon\zeta_t, \qquad t \in [0,T], \ \varepsilon \in \mathbb{R},$$

hence,

$$\zeta_t = \frac{d}{d\varepsilon} \mathbb{E}_{u+\varepsilon\zeta}[\xi_t - u_t]_{|\varepsilon=0}$$
$$= \frac{d}{d\varepsilon} \mathbb{E}[(\xi_t - u_t)\Lambda(u+\varepsilon\zeta)]_{|\varepsilon=0}$$
$$= \mathbb{E}\left[ (\xi_t - u_t) \frac{d}{d\varepsilon} \Lambda(u+\varepsilon\zeta)_{|\varepsilon=0} \right]$$
$$= \mathbb{E}_u\left[ (\xi_t - u_t) \frac{d}{d\varepsilon} \log \Lambda(u+\varepsilon\zeta)_{|\varepsilon=0} \right]$$

$$= \mathbb{E}_u\bigg[(\xi_t - u_t)\bigg(\int_0^T \frac{\dot{\zeta}_s}{\sigma_s^2}\,dX_s - \int_0^T \frac{\dot{\zeta}_s \dot{u}_s}{\sigma_s^2}\,ds\bigg)\bigg]$$

$$= \mathbb{E}_u\bigg[(\xi_t - u_t)\int_0^T \frac{\dot{\zeta}_s}{\sigma_s^2}\,dX_s^u\bigg]$$

$$= \mathbb{E}_u\bigg[(\xi_t - u_t)\int_0^t \frac{\dot{\zeta}_s}{\sigma_s^2}\,dX_s^u\bigg],$$

where the exchange between expectation and derivative is justified by classical uniform integrability arguments. Thus, by the Cauchy–Schwarz inequality and the Itô isometry, we have

$$\zeta_t^2 \leq \mathbb{E}_u\bigg[\bigg(\int_0^t \frac{\dot{\zeta}_s}{\sigma_s^2}\,dX_s^u\bigg)^2\bigg]\mathbb{E}_u[|\xi_t - u_t|^2]$$

$$= \int_0^t \frac{\dot{\zeta}_s^2}{\sigma_s^2}\,ds\,\mathbb{E}_u[|\xi_t - u_t|^2], \qquad t \in [0,T].$$

It then suffices to take

$$\zeta_t = \int_0^t \sigma_s^2\,ds, \qquad t \in [0,T],$$

to get

(5.4) $$\operatorname{Var}_u[\xi_t] = \mathbb{E}_u[|\xi_t - u_t|^2] \geq \int_0^t \sigma_s^2\,ds, \qquad t \in [0,T],$$

which leads to (3.5) after integration with respect to $\mu(dt)$. As noted above, $\hat{u} = (X_t)_{t\in[0,T]}$ is clearly unbiased under $\mathbb{P}_u$ and it attains the lower bound $\mathrm{R}(\sigma,\mu,\hat{u})$. $\square$

*Superefficient drift estimators.*

PROOF OF THEOREM 3.6.  Let $\delta\colon L^2(\Omega; H, \mathbb{P}_u) \to L^2(\Omega, \mathbb{P}_u)$ denote the closable adjoint of $D$, that is, the divergence operator under $\mathbb{P}_u$, which satisfies the integration by parts formula

(5.5)   $$\mathbb{E}_u[F\delta(v)] = \mathbb{E}_u[\langle v, DF\rangle], \qquad F \in \mathrm{Dom}(F),\ v \in \mathrm{Dom}(\delta),$$

with the relation

$$\delta(hF) = FX(h) - \langle h, DF\rangle;$$

see [7], for $F \in \mathrm{Dom}(D)$ and $h \in H$ such that $hF \in \mathrm{Dom}(\delta)$. Note that (5.5) is an infinite-dimensional version of the integration by parts (1.1), which can be proved using, for example, the countable sequence of Gaussian random variables appearing in the Paley–Wiener expansion (4.1) of $X$.

(i) We have
$$\mathbb{E}_u[\|X+\xi-u\|^2_{L^2([0,T],d\mu)}]$$
$$= \mathbb{E}_u\left[\int_0^T |X_t^u+\xi_t|^2 \mu(dt)\right]$$
$$= \mathbb{E}_u\left[\int_0^T |X_t^u|^2 \mu(dt)\right] + \|\xi\|^2_{L^2(\Omega\times[0,T],\mathbb{P}_u\otimes\mu)} + 2\mathbb{E}_u\left[\int_0^T X_t^u \xi_t \mu(dt)\right]$$
$$= \mathrm{R}(\gamma,\mu,\hat{u}) + \|\xi\|^2_{L^2(\Omega\times[0,T],\mathbb{P}_u\otimes\mu)} + 2\mathbb{E}_u\left[\int_0^T X_t^u \xi_t \mu(dt)\right]$$

and
$$\mathbb{E}_u[\xi_t X_t^u] = \mathbb{E}_u[\xi_t X^u(\chi_t)]$$
$$= \mathbb{E}_u[\xi_t \delta(\chi_t)]$$
$$= \mathbb{E}_u[\langle \chi_t, D\xi_t\rangle]$$
$$= \mathbb{E}_u[D_t \xi_t], \qquad t \in [0,T],$$

hence, (3.8) holds.

(ii) Next, we specialize the above argument to processes $\xi$ of the form
$$\xi_t = D_t \log F, \qquad t \in [0,T],$$
where $F$ is an a.s. strictly positive and sufficiently smooth random variable. From (3.8) we have
$$\mathbb{E}_u[\|X + D\log F - u\|^2_{L^2([0,T],\mu)}]$$
$$= \mathrm{R}(\gamma,\mu,\hat{u}) + \|D\log F\|^2_{L^2(\Omega\times[0,T],\mathbb{P}_u\otimes\mu)} + 2\mathbb{E}_u\left[\int_0^T D_t D_t \log F \mu(dt)\right]$$
$$= \mathrm{R}(\gamma,\mu,\hat{u}) + \mathbb{E}_u\left[\int_0^T \left(\left|\frac{D_t F}{F}\right|^2 + 2 D_t D_t \log F\right)\mu(dt)\right],$$

and we use the relation
$$\left|\frac{D_t F}{F}\right|^2 + 2 D_t D_t \log F = 2\frac{D_t D_t F}{F} - \left|\frac{D_t F}{F}\right|^2, \qquad t \in [0,T].$$

Concerning the last equality, we note that
$$2\frac{D_t D_t F}{F} - \left|\frac{D_t F}{F}\right|^2 = \frac{2}{\sqrt{F}} D_t\left(\frac{D_t F}{\sqrt{F}}\right) = \frac{4}{\sqrt{F}} D_t D_t \sqrt{F}, \qquad t \in [0,T],$$

which implies
$$(5.6) \qquad 4\frac{\Delta\sqrt{F}}{\sqrt{F}} = 2\frac{\Delta F}{F} - \int_0^T |D_t \log F|^2 \mu(dt),$$

and allows us to conclude from (i). $\square$

**Acknowledgments.** We thank the editors and referees for suggestions which led to several improvements of this paper, and, in particular, for communicating to us the references [2] and [11].

## REFERENCES

[1] ALÒS, E., MAZET, O. and NUALART, D. (2001). Stochastic calculus with respect to Gaussian processes. *Ann. Probab.* **29** 766–801. MR1849177

[2] BERGER, J. and WOLPERT, R. (1983). Estimating the mean function of a Gaussian process and the Stein effect. *J. Multivariate Anal.* **13** 401–424. MR0716932

[3] FOURDRINIER, D., STRAWDERMAN, W. E. and WELLS, M. T. (1998). On the construction of Bayes minimax estimators. *Ann. Statist.* **26** 660–671. MR1626063

[4] IBRAGIMOV, I. A. and ROZANOV, Y. A. (1978). *Gaussian Random Processes.* Springer, New York. MR0543837

[5] JAMES, W. and STEIN, C. (1961). Estimation with quadratic loss. *Proc. 4th Berkeley Sympos. Math. Statist. Probab.* **I** 361–379. Univ. California Press, Berkeley. MR0133191

[6] LIPTSER, R. S. and SHIRYAEV, A. N. (2001). *Statistics of Random Processes*, II. Springer, Berlin. MR1800858

[7] NUALART, D. (2006). *The Malliavin Calculus and Related Topics*, 2nd ed. Springer, Berlin. MR2200233

[8] PRIVAULT, N. and RÉVEILLAC, A. (2008). Stochastic analysis on Gaussian space applied to drift estimation. Preprint. arXiv:0805.2002v1.

[9] PRAKASA RAO, B. L. S. (1999). *Statistical Inference for Diffusion Type Processes.* Edward Arnold, London. MR1717690

[10] STEIN, C. (1981). Estimation of the mean of a multivariate normal distribution. *Ann. Statist.* **9** 1135–1151. MR0630098

[11] WOLPERT, R. and BERGER, J. (1982). Incorporating prior information in minimax estimation of the mean of a Gaussian process. In *Statistical Decision Theory and Related Topics III* **2** (*West Lafayette, Ind., 1981*) 451–464. Academic Press, New York. MR0705329

DEPARTMENT OF MATHEMATICS
CITY UNIVERSITY OF HONG KONG
83 TAT CHEE AVENUE
KOWLOON
HONG KONG
E-MAIL: nprivaul@cityu.edu.hk

LABORATOIRE DE MATHÉMATIQUES
UNIVERSITÉ DE LA ROCHELLE
AVENUE MICHEL CRÉPEAU
17042 LA ROCHELLE
FRANCE
E-MAIL: areveill@univ-lr.fr